\documentclass{amsart}
\usepackage{graphicx, palatino}

\def\cf{{\em cf. }}
\def\inv{^{-1}}
\def\Vt{\tilde{V}}
\def\C{\mathbb{C}}

\newtheorem{theorem}{Theorem}[section]
\newtheorem{proposition}[theorem]{Proposition}

\theoremstyle{definition}
\newtheorem{definition}[theorem]{Definition}

\theoremstyle{remark}
\newtheorem{remark}[theorem]{Remark}

\begin{document}

\title{A Remark on the Hard Lefschetz Theorem for K\"ahler Orbifolds}

\author{Wang Z. Z.}
\address{Fudan university\\ 
Institute of Mathematics\\ 
Shanghai 200433\\China}
\email{youxiang163wang@163.com}

\author{D. Zaffran}
\address{Fudan university\\ 
Institute of Mathematics\\ 
Shanghai 200433\\China}
\email{zaffran@fudan.edu.cn}

\subjclass[2000]{Primary 14F25; Secondary 53C12}

\commby{Jon Wolfson}

\begin{abstract}
We give a proof of the Hard Lefschetz Theorem for orbifolds that 
does not involve intersection homology. We use a foliated version of the Hard Lefschetz Theorem due to El Kacimi. 
\end{abstract}

\maketitle


\section*{Introduction} \noindent
Even though orbifolds are not smooth manifolds in general, these 
mildly singular objects retain a strong flavor of smoothness. 
For example, for a compact orientable orbifold $V$, 
Poincar\'e duality holds over the reals. If $V$ is moreover 
K\"ahler, it satisfies the Hard Lefschetz Theorem (HLT). 

The former was proved by Satake, by adapting the classical proof. 
However, the only known proof of the HLT for orbifolds uses the 
far-reaching version of the HLT in intersection homology. 
Fulton raises the problem of finding a more direct proof (\cf \cite{Ful}, p.105). 
We propose a solution based on a result of El Kacimi (\cf \cite{EK1}). 
We emphasize that the methods 
that we are using 
here are well known to the specialists of Riemannian foliations 
(\cf \cite{EK1}, \cite{EK2}, \cite{GHS}, \cite{Mo}, 
\cite{Rum}, \cite{Ton}). 

\section{The Resolution of an Orbifold as a Seifert bundle} \noindent
\begin{definition} \label{defKahlerOrbifolds}
Let $|V|$ be a
Hausdorff topological space. A {\em K\"{a}hler orbifold} structure $V$ on
$|V|$ is given by the following:

\noindent (i) An open cover ${\{V_i\}}_{i}$ of $|V|$.

\noindent (ii) For each $i\in I$, 
a connected and open $U_i\subset \mathbb{C}^n$ with a K\"{a}hler metric $h_i$; 
a finite subgroup $\Gamma_i$ of holomorphic isometries of $U_i$;
a continuous map $q_i:U_i\rightarrow V_i$, 
called a {\em local uniformization}, 
inducing a homeomorphism from
$\Gamma_i\setminus U_i$ onto $V_i$.

\noindent (iii) For all $x_i\in U_i$ and $x_j\in U_j$ such that
$q_i(x_i)=q_j(x_j)$, there exist 
$W_i \subset U_i$ and $W_j \subset U_j$, 
open connected neighbourhoods of $x_i$ and $x_j$, 
and a holomorphic isometry $\phi_{ji}: W_i \rightarrow W_j$, 
called a {\em change of charts}, 
such that $q_j \phi_{ji}=q_i$ on
$W_i$. 
\end{definition}
\noindent Remark: it is true, but not trivial, that 
$V$ is K\"ahler in the above ``orbifold''
sense if and only if, as a complex space, 
$V$ is K\"ahler in the sense of Grauert 
(defined in \cite{G}). 

\smallskip 
Let $V$ be a compact K\"{a}hler orbifold.
The representation of $V$ as the leaf space of a smooth foliation 
is due to \cite{GHS}, Prop. 4.1. We explain their construction, and check 
that the foliation is transversely K\"{a}hler. 

Fix $i$. For any 
local uniformization $q_i:U_i\rightarrow V_i$, 
we denote by $\pi_i: \tilde{U}_i \rightarrow U_i$ 
the principal bundle of unitary frames. 
The $\Gamma_i$-action on $U_i$ lifts naturally to $\tilde{U}_i$. 
Indeed, let $\xi_i=(\xi_i^1 ,\dots, \xi_i^n) \in \tilde{U}_i$ 
be a unitary frame at $x_i$ and let $g_i\in \Gamma_i$. 
Then $g_i\xi_i := \big( (g_i)_* \xi_i^1 ,\dots, (g_i)_*\xi_i^n \big)$ 
is a unitary frame at $g_i x_i$. As an isometry is determined 
by its derivative at one point, this lifted action is {\em free}. 
The quotient 
$\tilde{V}_i:=\Gamma_i \setminus \tilde{U}_i$ is a 
{\em $\mathcal{C}^{\infty}$-smooth manifold}, 
and  $\tilde{q}_i:\tilde{U}_i\rightarrow \tilde{V}_i$ is a finite 
non-ramified covering. 
 
\begin{figure}
\includegraphics{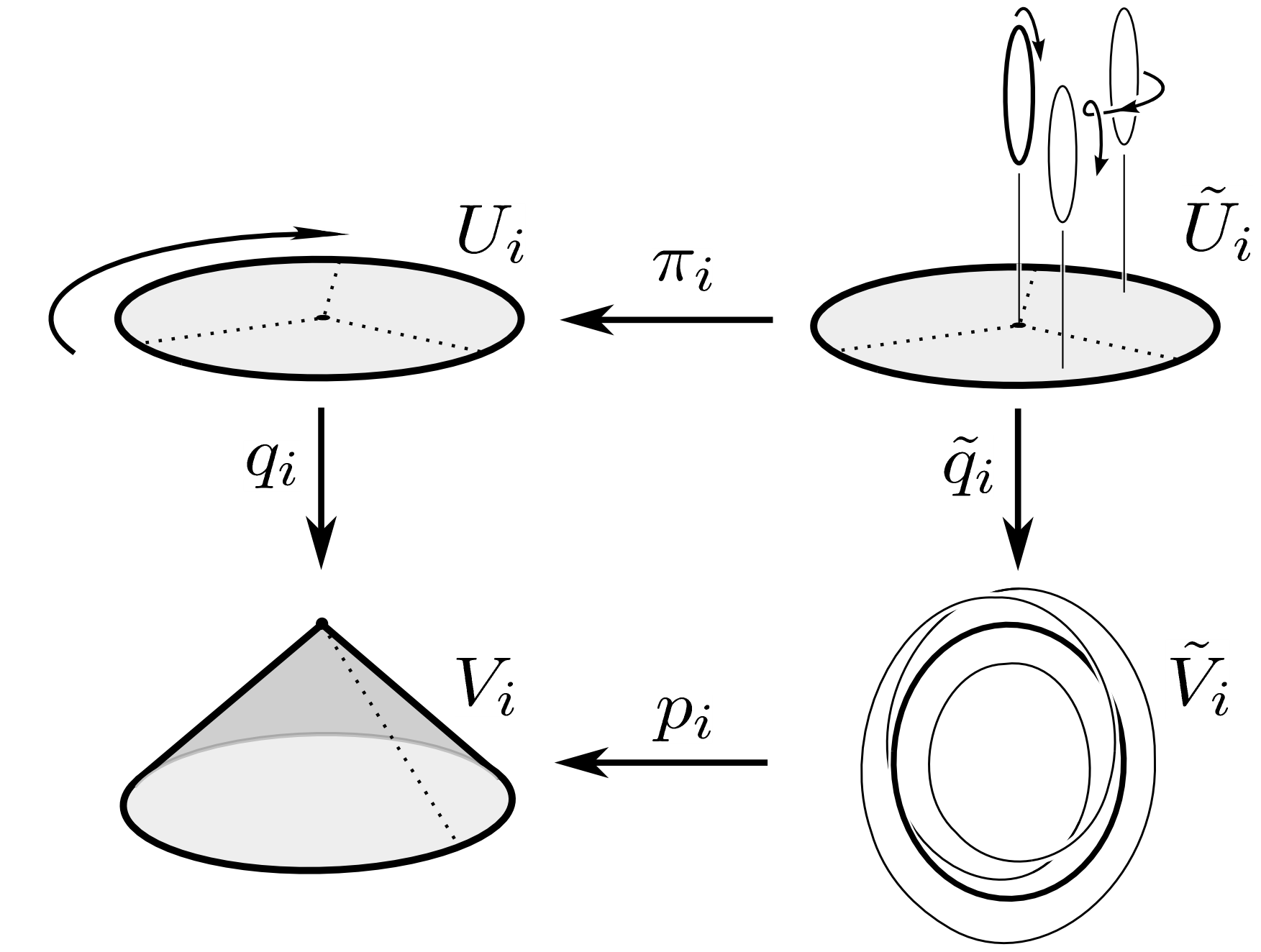}
\end{figure}

Fix $i$ and $j$. Any change of charts 
$\phi_{ji}: W_i \subset U_i \rightarrow W_j \subset U_j$ 
lifts, as above, to a diffeomorphism 
$\tilde{\phi}_{ji}: 
\tilde{W}_i\subset \tilde{U}_i  \rightarrow \tilde{W}_j \subset \tilde{U}_j, 
\ \xi_i \mapsto \big( 
({\phi}_{ji})_* \xi_i^1 ,\dots, ({\phi}_{ji})_*\xi_i^n \big)$. 
\begin{remark} 
For any two changes of charts $\phi_{ji}$ and $\phi'_{ji}$ defined 
on the same $W_i \subset U_i$, 
there exists $g_j \in \Gamma_j$ such that $\phi'_{ji}= g_j \phi_{ji}$ 
(this follows from the classical case of a non-ramified $q_j$). 
Accordingly, 
$\tilde{\phi}'_{ji}= g_j \tilde{\phi}_{ji} 
\text{ on } \tilde{W}_i \subset \tilde{U}_i.$
\label{equiv} 
\end{remark} 
\noindent Denote $V_{ji}:= V_i \cap V_j$. 
We want to define a lifting of the obvious gluing map 
$$f_{ji}: V_{ji} \subset {V}_i \rightarrow V_{ji} \subset {V}_j$$ 
to a gluing diffeomorphism 
$$\tilde{f}_{ji}: 
p_i\inv(V_{ji}) \subset \tilde{V}_i \rightarrow 
p_j\inv(V_{ji}) \subset \tilde{V}_j.$$
Let $z_i\in p_i\inv(V_{ji})$. 
Let $\phi_{ji}: W_i \rightarrow W_j$ 
be a change of charts with $y_i\in \tilde{W}_i$ such that 
$\tilde{q}_i(y_i) = z_i$. 
Take $\tilde{f}_{ji}(z_i):= \tilde{q}_j \tilde{\phi}_{ji} (y_i)$. 
This is well-defined: if we choose another point $g_i y_i$ and 
another change of charts 
$\phi'_{ji}: g_iW_i \rightarrow W'_j$, then both 
$\phi_{ji}$ and $\phi'_{ji} g_i$ are changes of charts on ${W}_i$, 
so, by {\em Rem.} \ref{equiv}, there is a $g_j\in\Gamma_j$ such that 
$\tilde{\phi}'_{ji} g_i = g_j \tilde{\phi}_{ji}$. Thus 
$\tilde{q}_j \tilde{\phi}'_{ji} g_i (y_i) = 
 \tilde{q}_j g_j \tilde{\phi}_{ji} (y_i) =  
 \tilde{q}_j \tilde{\phi}_{ji} (y_i)$.

Also from {\em Rem.} \ref{equiv}, for all 
$i,j,k,\ \tilde{f}_{ki}=\tilde{f}_{kj} \tilde{f}_{ji}$. 
Therefore, the disjoint collection $\{\tilde{V_i}{\}}_i$ 
glued according to $\{\tilde{f}_{ji}{\}}_{ij}$ yields a smooth 
compact manifold that we 
call $\tilde{V}$. It is Hausdorff as a consequence of 
the following
\begin{proposition} For any compact orbifold $V$, the manifold $\tilde{V}$ 
constructed as above admits a $U(n)$-action, turning 
$p: \tilde{V} \rightarrow V$ into a Seifert principal bundle. 
\end{proposition}
\begin{proof}
First we define, for any fixed $i$, a right $U(n)$-action on $\tilde{U}_i$. 
Let $\xi_i=(\xi_i^1 ,\dots, \xi_i^n) \in \tilde{U}_i$, and 
$A=(A^{\alpha \beta})_{1\leq \alpha ,\beta \leq n}\in U(n)$. 
We define $\xi_i A$ as the ``matrix product'' 
$(\xi_i^1 ,\dots, \xi_i^n) A$. 
This means that 
$\xi_i A:=(\zeta^1, \dots, \zeta^n)$ where for all $\beta$, 
$\zeta^{\beta}:= 
{
\sum_{\alpha}} A^{\alpha\beta} \xi_i^{\alpha}$.
Now, by linearity of the derivative, 
for all 
$g_i, A$ and $\xi_i,\quad g_i (\xi_i A) = (g_i \xi_i) A$. 
So the action descends to $\tilde{V_i}$. Similarly, 
the $U(n)$-actions on $\tilde{V_i}$ and $\tilde{V_j}$ 
commute with $\tilde{f}_{ji}$, endowing $\tilde{V}$ with 
a global action. 

The Seifert bundle structure follows from the construction. 
\end{proof}
\begin{proposition} \label{tK}The foliation $\mathfrak{F}$ 
of $\tilde{V}$ by $U(n)$-orbits is transversely K\"{a}hler, 
i.e., 
$\mathfrak{F}$ is transversely holomorphic (so 
the normal bundle admits a complex structure $J$), 
and there exists 
a closed real 2-form $\omega$ on $\tilde{V}$ such that:

(i) the tangent bundle of $\mathfrak{F}$ is the kernel of $\omega$;

(ii) The quadratic form
$h(-,-):=\omega(J-,-)+\sqrt{-1}\omega(-,-)$ defines a
Hermitian metric on the normal bundle of $\mathfrak{F}$. 
\end{proposition}
\begin{proof}
For an open subset $O \subset \Vt$ small enough, we can assume 
that $O \subset \tilde{q}_i(\tilde{V}_i)$ for some $i$, 
and take a smooth section $\sigma: O \rightarrow \tilde{V}_i$ of $\tilde{q}_i$. 
Then, on $O$, the foliation $\mathfrak{F}$ is defined by the submersion 
$\pi_i\sigma : O \rightarrow U_i$, which induces an isomorphism between 
the normal bundle of $\mathfrak{F}$ and the tangent bundle of $U_i$. 
By hypothesis, we have a $\Gamma_i$-invariant K\"ahler form $\omega_i$ 
on $U_i$. Define $\omega:=\pi_i^{*}\omega_i$, a well-defined form on $\Vt$ that satisfies the above properties.
\end{proof}
\section{A taut and transversely K\"ahler metric}

\noindent
Let $g_0$ be a Riemannian metric on $\tilde{V}$. 
Up to averaging it over the action, 
we can assume that it is $U(n)$-invariant. 
Denote $m= \dim U(n)$. Fix a global frame on $U(n)$ given by right-invariant vector fields $E_1,\dots,E_m$.
Let $e_1, \dots, e_m$ be the associated 
fundamental vector fields on $\Vt$, and 
for any $[\xi] \in \Vt$, let 
$M_0[\xi]:= \Big( g_0\big(e_k[\xi],e_l[\xi]\big) 
{\Big)}_{1\leq k,l \leq m} \in \mathbb{R}^{m\times m}$, 
and define a smooth function $u$ on $\Vt$ by 
$$
u_0[\xi]:=\big(\det M_0[\xi]\big)^{-\frac{1}{m}}.
$$
By invariance, $M_0$ and $u_0$ are constant on each orbit. 
Now, we define a new Riemannian metric on $\Vt$ by the conformal change
$g_1=u_0 g_0$.
\begin{proposition} \label{taut}
With respect to $g_1$, the $U(n)$-orbits are minimal submanifolds.
\end{proposition}
\begin{proof} (inspired by the proof of \cite{Ton} Prop 7.6.) 
Since this is a local statement, we can assume that 
we are on $\Vt_i$ for some $i$. 
Using the covering map $\tilde{q}_i$, 
we lift the metric $g_1$ to a metric $h_1$ on $\tilde{U}_i$. 
Let 
 $X_0 \subset \tilde{U}_i$ be any $U(n)$-orbit, 
and let $Z$ be an arbitrary vector field 
that is compactly supported and normal to $X_0$. 
Let ${ \{ \varphi_t \} }_{t\in \mathbb{R}}$ be the 
flow associated to $Z$, and define 
$X_t:= \varphi_t(X_0)$. 
By the first variational formula (or \cite{Jos} 3.6.4), 
it is enough to prove that 
$$ 
\frac{d}{dt}
\big{(\textrm{Vol}\ X_t \big)}_{|_{t=0}}=0
.
$$
Because the metric $h_1$ is $U(n)$-invariant, we can replace $Z$ with its average 
over the action. This new $Z$ is still orthogonal to $X_0$, 
and is $U(n)$-invariant, so $X_t$ is a $U(n)$-orbit for every $t$. 
On $\Vt$, the volume form of 
a $U(n)$-orbit is 
$$\big( \det M_1 \big)^{\frac{1}{2}}\ \check{e}_1 \wedge\dots\wedge \check{e}_m,$$ 
where $M_1= \Big( g_1\big(e_k,e_l\big) 
{\Big)}_{1\leq k,l \leq m}$, and 
$\check{e}_1, \dots, \check{e}_m$ is the dual basis. 
On the other hand, 
$\big( \det M_1 \big)^{\frac{1}{2}}
= u_0^{\frac{m}{2}} \big( \det M_0 \big)^{\frac{1}{2}}
= \big( (\det M_0)^{-\frac{1}{m}} \big)^{\frac{m}{2}} \big( \det M_0 \big)^{\frac{1}{2}}
=1.$

Fix an arbitrary $t$. Choosing any 
point on $M_t$, we get a one-to-one a parameterization 
$f: U(n) \rightarrow X_t$. 
Then 
$$\textrm{Vol}\ X_t 
=\int_{U(n)} f^*\tilde{q}_i^*\Big(\check{e}_1 \wedge\dots\wedge \check{e}_m\Big)
=\int_{U(n)} \check{E}_1\wedge\dots\wedge\check{E}_m,$$
which does not depend on $t$.
\end{proof}

The metric $g_1$ determines an embedding 
$N\subset T\Vt$ of the normal bundle of $\mathfrak{F}$. 
Modifying $g_1$, we define a metric $g$ on $\Vt$ by replacing 
${g_1}_{\vert N}$ with the metric $h$ obtained in Prop. \ref{tK}. 
The proof of Prop. \ref{taut} shows that the leaves 
of $\mathfrak{F}$ are still minimal submanifolds with respect to $g$. 

\section{Main result}

\noindent
Let $M$ be any compact smooth manifold, with a foliation $\mathfrak{F}$. 
A form $\alpha \in \Omega^*(M,\C)$ is called {\em basic} when for all 
vectors $Z$ tangent to the foliation, $\iota_Z \alpha = 0 = \iota_Z d\alpha$, 
where $\iota_Z$ denotes the interior product. 
The basic forms determine a subcomplex of the usual de Rham 
complex of $M$. The associated cohomology groups are called the 
basic cohomology groups, and are denoted $H_B^*(\mathfrak{F},\C)$. 
They can be thought of as 
a substitute for the de Rham cohomology of the leaf space, which can be 
very singular and/or not Hausdorff in general. 
When the leaf space is an orbifold $V$, the basic cohomology is 
isomorphic to the singular cohomology of $V$ over $\C$ 
\mbox{\text{(\cf \cite{Pfl} 5.3)}}. 

Let $V$ be a compact K\"ahler orbifold, with $\Vt$, 
$\mathfrak{F}$ and $g$ as described in the previous sections. 
We want apply El Kacimi's Basic Hard Lefschetz Theorem: 
if $\mathfrak{F}$ is transversely K\"ahler with respect to $\omega$, 
and $H_B^{2n}(\mathfrak{F},\C)\neq 0$, 
then for any $k$, $\alpha \mapsto \alpha \wedge \omega^k$ is 
an isomorphism 
$H_B^{n-k}(\mathfrak{F},\C) \rightarrow H_B^{n+k}(\mathfrak{F},\C)$, 
where $n$ is 
the complex codimension of $\mathfrak{F}$. We refer to \cite{EK1} 
for a more complete statement of the theorem. 

By a theorem of Masa (\cf \cite{Ma}), the existence of a metric 
turning the leaves 
into minimal submanifolds is equivalent to $H_B^{2n}(\mathfrak{F},\C)\neq 0$. 
This completes the proof of the Hard Lefschetz Theorem for $V$. 

Remark: in order to avoid 
the use of \cite{Ma}, one can easily 
adapt the proof of \cite{EK1} (written when Masa's theorem was only a conjecture). 
The assumption $H_B^{2n}(\mathfrak{F},\C)\neq 0$ 
is used in \cite{EK1} to prove the classical Hodge-theoretic adjoint formula: 
$\delta= (-1)^i *\partial*$. Due to the vanishing of the mean curvature 
associated to our metric $g$, this formula follows 
from \cite{Ton} Th. 7.10.


\bigskip


\begin{thebibliography}{999999}

\bibitem [GHS]{GHS} J.Girbau, A.Haefliger, D.Sundararaman, 
{\it On deformations of transversely holomorphic foliations}, 
J. Reine Angew. Math. {\bf 345} (1983), 122--147. 

\bibitem [G]{G} H. Grauert, 
{\it \"Uber Modifikationen und exzeptionelle analytische Mengen}, 
Math. Ann. {\bf 146} (1962), 331--368.

\bibitem [EK1]{EK1} A. El Kacimi Alaoui, 
{\it Op\'{e}rateurs transversalement elliptiques sur un 
feuilletage riemannien et applications}, 
Compositio Mathematica {\bf 73} (1990), 57--106.

\bibitem [EK2]{EK2} A. El Kacimi-Alaoui, 
{\it Stabilit\'e des $V$-vari\'et\'es kahl\'eriennes}, 
Lecture Notes in Math. {\bf 1345}, 111--123. Springer, 1988.

\bibitem [Ful]{Ful} W. Fulton, {\it Introduction to toric varieties}, 
Princeton University Press, 1993. 

\bibitem [Jos]{Jos} J. Jost, {\it Riemannian geometry and geometric analysis}, 
Third edition. Universitext. Springer, 2002.

\bibitem [Ma]{Ma} X. Masa, 
{\it Duality and minimality in Riemannian foliations}, 
Commentarii Mathematici Helvetici {\bf 67} (1992), 17--27.

\bibitem[Mo]{Mo} P. Molino, {\it Riemannian foliations}, 
Progress in Mathematics {\bf 73}. Birkhäuser, Boston, 1988.

\bibitem[Pfl]{Pfl} M. Pflaum, 
{\it Analytic and geometric study of stratified spaces}, 
Lecture Notes in Mathematics {\bf 1768}. Springer, 2001.

\bibitem[Rum]{Rum} H. Rummler, 
{\it Quelques notions simples en g\'eom\'etrie riemannienne et leurs applications aux feuilletages compacts}, 
Comment. Math. Helv. {\bf 54} (1979), 224--239. 

\bibitem [Ton]{Ton} P. Tondeur,
{\it Geometry of foliations}, 
Monographs in Mathematics {\bf 90}. Birkh\"auser Verlag, Basel, 1997.


\end{thebibliography}
\end{document}